\documentclass{article}
\usepackage{CJK,CJKnumb,CJKulem,times,dsfont,ifthen,mathrsfs,latexsym,amsfonts,color}
\usepackage{amsmath,amsthm,makeidx,fontenc,amssymb,bm,graphicx,psfrag,listings,curves,extarrows}
\makeindex
\newtheorem{definition}{Definition}[section]
\newtheorem{theorem}{Theorem}[section]
\newtheorem{lemma}{Lemma}[section]
\newtheorem{corollary}{Corollary}[section]
\newtheorem{proposition}{Proposition}[section]
\newtheorem{remark}{Remark}[section]

\newcommand{\s}{\section}

\newcommand{\R}{\mathbb R}

\newcommand{\lab}{\label}
\newcommand{\bt}{\begin{theorem}}
\newcommand{\et}{\end{theorem}}
\newcommand{\bl}{\begin{lemma}}
\newcommand{\el}{\end{lemma}}
\newcommand{\bd}{\begin{definition}}
\newcommand{\ed}{\end{definition}}
\newcommand{\bc}{\begin{corollary}}
\newcommand{\ec}{\end{corollary}}
\newcommand{\bp}{\begin{proof}}
\newcommand{\ep}{\end{proof}}
\newcommand{\bx}{\begin{example}}
\newcommand{\ex}{\end{example}}
\newcommand{\bi}{\begin{exercise}}
\newcommand{\ei}{\end{exercise}}
\newcommand{\bo}{\begin{proposition}}
\newcommand{\eo}{\end{proposition}}
\newcommand{\br}{\begin{remark}}
\newcommand{\er}{\end{remark}}
\newcommand{\be}{\begin{equation}}
\newcommand{\ee}{\end{equation}}
\newcommand{\ba}{\begin{align}}
\newcommand{\ea}{\end{align}}
\newcommand{\bn}{\begin{enumerate}}
\newcommand{\en}{\end{enumerate}}
\newcommand{\bg}{\begin{align*}}
\newcommand{\bcs}{\begin{cases}}
\newcommand{\ecs}{\end{cases}}

\newcommand{\bean}{\begin{eqnarray*}}
\newcommand{\eean}{\end{eqnarray*}}

\numberwithin{equation}{section}

\begin{document}
\begin{CJK*}{GBK}{song}
\title{\bf  Existence of extremal functions  for  a family of Caffarelli-Kohn-Nirenberg inequalities\thanks{Supported by NSFC (11025106, 11371212, 11271386) and the Both-Side Tsinghua Fund.
E-mails: zhongxuexiu1989@163.com\quad \quad wzou@math.tsinghua.edu.cn}}

\date{}
\author{
{\bf  X. Zhong\;     \&  \;  W. Zou}\\
\footnotesize Department of Mathematical Sciences, Tsinghua University,  Beijing 100084, China}

\maketitle

\vskip0.6in

\begin{center}
\begin{minipage}{120mm}
\begin{center}{\bf Abstract}\end{center}

Consider  the following inequalities due to Caffarelli, Kohn  and Nirenberg {\it (Compositio Mathematica,1984):}
$$\Big(\int_\Omega \frac{|u|^r}{|x|^s}dx\Big)^{\frac{1}{r}}\leq C(p,q,r,\mu,\sigma,s)\Big(\int_\Omega \frac{|\nabla u|^p}{|x|^\mu}dx\Big)^{\frac{a}{p}}\Big(\int_\Omega \frac{|u|^q}{|x|^\sigma}dx\Big)^{\frac{1-a}{q}},$$
where $\Omega \subset \R^N (N\geq 2)$ is an open set; $p, q,  r, \mu, \sigma, s, a$ are some parameters satisfying some balanced  conditions.  When $\Omega$ is a cone in $\R^N$ (for example, $\Omega=\R^N)$, we   prove the sharp constant $C(p,q,r,\mu,\sigma,s)$ can be achieved  for a very large parameter space.  Besides, we  find some sufficient conditions which guarantee  that the following  Sobolev spaces
$$W_{\mu}^{1,p}(\Omega),\;  W_{\mu}^{1,p}(\Omega)\cap L^p(\Omega), \; H^{1,p}(\R^N) $$  are   compactly embedded into $L^r(\R^N, \frac{dx}{|x|^s})$ for some new ranges  of parameters, where
$\displaystyle W_{\mu}^{1,p}(\Omega)$ is the completion of $C_0^\infty(\Omega)$ with respect to the norm
$\displaystyle \Big(\int_\Omega \frac{ |\nabla u|^p}{|x|^\mu}dx\Big)^{\frac{1}{p}}. $
As applications, we also study the equation  $$\displaystyle -div\Big(\frac{|\nabla u|^{p-2}\nabla u}{|x|^\mu}\Big)=\lambda V(x)|u|^{q-2}u, \;\;\;   u\in W_{\mu}^{1,p}(\Omega)$$  under some proper conditions on $V(x)$.


 \vskip0.023in

{\it   Key  words:}   Rellich-Kondrachov theorem, Caffarelli-Kohn-Nirenberg inequality,   Ground state,   Extremal functions.

\vskip0.023in

\end{minipage}
\end{center}
\vskip0.26in
\newpage
\s{Introduction}
\renewcommand{\theequation}{1.\arabic{equation}}
\renewcommand{\theremark}{1.\arabic{remark}}
\renewcommand{\thedefinition}{1.\arabic{definition}}
In 1984,  Caffarelli,   Kohn and   Nirenberg established  a family of interpolation inequalities,  nowadays  called Caffarelli-Kohn-Nirenberg  (CKN) inequalities.
\vskip 0.03in
\noindent
{\bf Theorem A}\lab{2013-9-10-th1}(cf.\cite{CaffarelliKohnNirenberg.1984}) {\it
Assume that $p, q, r,  \alpha, \beta, \sigma $ and $a$ are fixed real numbers (called parameters) satisfying
\be\lab{2013-9-10-e1}
p\geq 1,   \;\;  q\geq 1, \;\; r>0,  \;\; 0\leq a\leq 1;
\ee
\be\lab{2013-9-10-e2}
\frac{1}{p}+\frac{\alpha}{N}>0,\;\;  \frac{1}{q}+\frac{\beta}{N}>0,\;\; \frac{1}{r}+\frac{\gamma}{N}>0,
\ee
where $
\gamma=a\sigma+(1-a)\beta. $
Then there exists a positive constant $C$ such that the following inequality holds
\be\lab{2013-9-10-e4}
\big||x|^\gamma u\big|_{L^r}\leq C\big||x|^\alpha|Du|\big|_{L^p}^{a}\big||x|^\beta u\big|_{L^q}^{1-a}, \quad \forall  u\in C_0^\infty(\R^N)
\ee
if and only if the following relations hold:
\be\lab{2013-9-10-e5}
\frac{1}{r}+\frac{\gamma}{N}=a\big(\frac{1}{p}+\frac{\alpha-1}{N}\big)+(1-a)\big(\frac{1}{q}+\frac{\beta}{N}\big)
\ee
(this is dimensional balance)
$$0\leq \alpha-\sigma\quad \hbox{if}\;a>0,$$
and
$$\alpha-\sigma\leq 1\quad \hbox{if}\;a>0\;\hbox{and}\;\frac{1}{p}+\frac{\alpha-1}{N}=\frac{1}{r}+\frac{\gamma}{N}.$$
Furthermore, on any compact set in parameter space in which (\ref{2013-9-10-e1}),(\ref{2013-9-10-e2}),(\ref{2013-9-10-e5}) and
$0\leq \alpha-\sigma\leq 1$ hold, the constant $C$ is bounded.}

\hfill$\Box$

Some  variant versions of  the CKN inequality with  higher order derivatives   were given by Lin  \cite{Lin.1986}.   Note that the CKN inequality and its variance include many well-known inequalities such as the Hardy-Sobolev inequality, Gagliardo-Nirenberg inequality, etc.   They play  a crucial  role  in the  elliptic partial differential equations.  Recall  a  version of the Gagliardo-Nirenberg inequality
 \be\lab{2014-9-5-e7}
    |u|_r\leq C|\nabla u|_2^a|u|_{2}^{1-a}.
    \ee
    When $2<r<2^*:=\frac{2N}{N-2}\;(N\geq 3)$, the dimensional balance condition implies that $0<a<1$. Then by the Young inequality and Sobolev inequality, we see that
    $\displaystyle H^1(\R^N)\hookrightarrow L^r(\R^N)$ is an continuous embedding for $2\leq r\leq 2^*$ which has been wildly used now.   If we consider that $a=1$ in (\ref{2013-9-10-e4}), then we have the following inequality without  interpolation:
\be\lab{2013-9-05-e5}
\big\||x|^\gamma u\big\|_{L^r}\leq C\big\||x|^\alpha |Du|\big\|_{L^p}\;\hbox{for all}\; u\in C_0^\infty(\R^N),
\ee
where $p\geq 1,  \frac{1}{r}+\frac{\gamma}{N}>0,  r\begin{cases}\leq p^*:=\frac{Np}{N-p},\quad &\hbox{if}\;N>p,\\ <\infty,&\hbox{if}\;N\leq p,\end{cases}$ and
\be\lab{2013-9-05-e6}
\frac{1}{r}+\frac{\gamma}{N}=\frac{1}{p}+\frac{\alpha-1}{N},
\ee
which is dimensional balance condition. We call (\ref{2013-9-05-e5})  the general Hardy-Sobolev inequality since when $\gamma=\alpha=0$,  (\ref{2013-9-05-e5})  returns to  the classical Sobolev inequality:
\be\lab{2013-9-05-e7}
|u|_r\leq C|\nabla u|_p\;\hbox{with}\;r=p^*, N>p,
\ee
When $\alpha=-t, \gamma=-t-1$, (\ref{2013-9-05-e5}) becomes
\be\lab{2013-9-05-e8}
\big|\frac{u}{|x|^{t+1}}\big|_p\leq C\big|\frac{\nabla u}{|x|^t}\big|_p,  \;\; N-p-pt>0,
\ee
which is called  the general weighted Hardy inequality.

\vskip0.1in
  Much progress has been made on (\ref{2013-9-05-e5})   for the case of $p=2$. For example,
in \cite{Aubin.1976,Talenti.1976a},   Aubin and Talenti gave the best constant and the minimizers for the Sobolev inequality (\ref{2013-9-05-e7})   via   the Schwarz symmetrization and the Bliss inequality in \cite{Bliss.1930}.   In \cite{Lieb.1983}, Lieb  applied  the same type of symmetrization to study (\ref{2013-9-05-e5}) with $\alpha=0, p=2, -1<\gamma<0$.   The results of \cite{Lieb.1983}  had been    generalized by Chou and Chu
in \cite{ChouChu.1993} to the case of $\alpha-1<\gamma\leq \alpha\leq 0, p=2$. A  further generalization was  given by Wang and Willem in  \cite{WangWillem.2000} for the case of $p=2$. When $p=2$ and $\alpha>0$, it  was  also studied in the papers \cite{CatrinaWang.2001}. For the case of $p\neq 2$ but with different geometries of the domain $\Omega\subset \R^N$,  we refer to \cite{BartschPengZhang.2007}. More  results about the related progress,  we refer to \cite{WangWillem.2003,LinWang.2004,DolbeaultEstebanLossEtAl.2009,DolbeaultEsteban.2012a,DolbeaultEsteban.2011}.  We remark that the   papers mentioned in this paragraph  mainly deal with the  inequality (\ref{2013-9-05-e5}) without  interpolation term.

\vskip0.1in

When $a\neq 1$, the CKN  inequality involves three terms (i.e., interpolation), which make  the problem much tough and there are rare paper investigating  this case, we just find the following
partial answers (see  the review paper by Dolbeault and   Esteban \cite{DolbeaultEsteban.2012b}):

 \begin{itemize}

\item When $\alpha=\beta=\gamma=0, p=2, q=p+1$ 
and $ r=2p$. For such a very special case, the sharp constant and the extremal functions of inequality \eqref{2013-9-10-e4} are given by  Del Pino and  Dolbeault \cite{DelDolbeault.2002}.

 \item When  $p=q=2,  -\frac{N-2}{2}<\alpha, \beta=\alpha-1, \alpha-1\leq \gamma<\alpha,$ and $ r=\frac{2N}{N+2(\gamma-\alpha)}$.   Under these assumptions,  together with a special region of $a$ and  other conditions,  the sharp constant and  extremal functions of  the CKN  inequality (\ref{2013-9-10-e4}) are studied by Dolbeault,  Esteban,    Tarantello  and  Tertikas \cite{DolbeaultEstebanTarantelloTertikas.2011},  Dolbeault  and Esteban \cite{DolbeaultEsteban.2012}.
\end{itemize}
 In the current paper,  we consider the general cases of  the CKN  inequality: $p>1$ and it has   interpolation term.

\vskip 0.13in

We   make a   transformation first.  Let $\alpha=-\frac{\mu}{p}, \beta=-\frac{\sigma}{q}, \gamma=-\frac{s}{r}$ in (\ref{2013-9-10-e4}), then a direct computation shows that $a=\frac{[(N-\sigma)r-(N-s)q]p}{[(N-\sigma)p-(N-\mu-p)q]r}$. We obtain  the following version of the CKN  inequality:
 \begin{align}\lab{2014-6-5-xe4}
 \Big(\int_\Omega \frac{|u|^r}{|x|^s}dx\Big)^{\frac{1}{r}}\leq C(p,q,r,\mu,\sigma,s)\Big(\int_\Omega \frac{|\nabla u|^p}{|x|^\mu}dx\Big)^{\frac{a}{p}}\Big(\int_\Omega \frac{|u|^q}{|x|^\sigma}dx\Big)^{\frac{1-a}{q}},
 \end{align}
In present paper,   when $\Omega$ is a cone (i.e., $\lambda x\in \Omega$ for all $x\in \Omega$ and $\lambda>0$), we can obtain the existence of extremal functions for the CKN  inequalities  (\ref{2014-6-5-xe4}).
Define $$p^*(s, \mu):=\frac{p(N-s)}{N-p-\mu}.$$


\bt\lab{2014-10-14-mainth1}
Let $\Omega\subset \R^N$ be a cone (in particular,  $\Omega=\R^N$).
Assume that $\displaystyle p>1, s>0,   \; \max\{\sigma, s\}< \mu+p<N$, $\displaystyle  1\leq r, 1\leq q<\min\{p^*,p^*(\sigma, \mu)\}$, $\displaystyle \max\{\frac{p(\sigma-s)}{N-\mu-p}+q, \frac{\sigma-s}{N-\sigma}q+q\}<r<\min\{p^*, p^*(s, \mu)\}$ and
\be\lab{2014-9-30-bu2}
\begin{cases}
p(s-\sigma)+q(\mu+p-s)< r(\mu+p-\sigma)\\
(Np-Nr+pr)(s-\sigma)> (N\mu-Ns+ps)(r-q)
\end{cases},
\ee
then the sharp constant $C(p,q,r,\mu,\sigma,s)$ can be achieved and
    $$C(p,q,r,\mu,\sigma,s)=\big(\frac{1}{\rho}\big)^{\frac{(\mu+p-\sigma)r+(p-q)(N-s)}{r[(N-\sigma)p-(N-\mu-p)q]}},$$
where
\be\lab{2014-9-21-we1}\rho:=\inf\Big\{\int_\Omega \frac{|\nabla u|^p}{|x|^\mu}dx+\lambda^*\int_\Omega \frac{|u|^q}{|x|^\sigma}dx\;:\; \int_{\Omega}\frac{|u|^r}{|x|^s}dx=1\Big\}\ee
which can be attained    and
            \begin{align}\lab{2014-9-21-we2}
             \lambda^*:=&\Big\{\frac{p(N-s)-(N-\mu-p)r}{(\mu+p-\sigma)r+(p-q)(N-s)}\Big\}^{\frac{(\mu+p-\sigma)r+(p-q)(N-s)}{(N-s)p-(N-\mu-p)r}}\nonumber\\
             &\cdot\Big\{\frac{(N-\sigma)r-(N-s)q}{p(N-s)-(N-\mu-p)r}\Big\}^{\frac{(N-\sigma)r-(N-s)q}{p(N-s)-(N-\mu-p)r}}.
             \end{align}
\et

\br\lab{2014-10-25-R1}  When  $\sigma=0, 1<p=q<N$, each of the  the following conditions meets
the hypotheses of Theorem   \ref{2014-10-14-mainth1}:
\begin{itemize}
\item[(1)]$\mu=0, 0<s<p<N, p< r<\frac{p(N-s)}{N-p}$;
\item[(2)]$\mu>0, \frac{N\mu(r-p)}{p^2}< s<\mu+p<N, p<r<\min\{\frac{pN}{N-p}, \frac{p(N-s)}{N-\mu-p}\}$;
\item[(3)] $\mu<0, 0<s<\mu+p<N, p< r< \frac{p(N-s)}{N-\mu-p}$.
\end{itemize}
In fact, under these conditions we shall  show that the embedding $$W_{\mu}^{1,p}(\Omega)\cap L^p(\Omega)\hookrightarrow L^r(\Omega, \frac{dx}{|x|^s})$$
is a compact embedding, where we denote by $W_{\mu}^{1,p}(\Omega)$ the completion of $C_0^\infty(\Omega)$ with respect to the norm
\be\lab{zz=z}\|u\|:= \big(\int_\Omega\frac{|\nabla u|^p}{|x|^\mu}dx\big)^{\frac{1}{p}}\ee
and $L^r(\Omega, \frac{dx}{|x|^s})$ stands for the completion of $C_0^\infty(\Omega)$ with respect to the norm of
$$|u|_{r, s}:=\big(\int_\Omega \frac{|u|^r}{|x|^s}dx\big)^{\frac{1}{r}}.$$
See Corollary \ref{2014-10-1-cro1} in Section 3.

\er

\br
   It is well known that $H^{1,p}(\R^N)\hookrightarrow L^r(\R^N)$ is a continuous embedding for $r\in [p, p^*)$ but not compact.  However,  we will prove  that $H^{1,p}(\R^N)\hookrightarrow L^r(\R^N, \frac{dx}{|x|^s})$ is a  compact embedding for $s>0$ and $ r\in [p, p^*(s))$.  See  Remark
   \ref{2014-9-19-r1}.
\er

\vskip 0.13in
In this paper, we also
study the following problem
\be\lab{2014-9-7-xe2}
-div\big(\frac{|\nabla u|^{p-2}\nabla u}{|x|^\mu}\big)=\lambda V(x)|u|^{q-2}u,\;\;  u\in W_{\mu}^{1,p}(\Omega),
\ee
where $1<q<p^*:=\frac{Np}{N-p}, 1<p<N, \mu+p<N$ and $\Omega \subset \R^N$ is an open   Lipschitz domain     and $\lambda$ is a parameter.
We assume that $V\in L_{loc}^1(\Omega), V=V_+-V_-$, where $\displaystyle V_\pm(x):=\max\{\pm V(x), 0\}$.
Our basic assumption is
\begin{itemize}
\item [${\bf (H)}$]  $V\in L_{loc}^1(\Omega), V_+=V_1+V_2\neq 0,$ there exists some $\frac{N\mu}{N-p}\leq \eta< \min\{\mu+p, N-\frac{q}{p}(N-\mu-p)\}$ such that $\big|V_1(x)\big|^{\frac{p^*(\eta, \mu)}{p^*(\eta,\mu)-q}}|x|^{\frac{q\eta}{p^*(\eta,\mu)-q}}\in L^1(\Omega)$ and one of the following holds
    \begin{itemize}
    \item[$(H_1)$] $1<q<p, \Omega$ is bounded and $\displaystyle \sup_{y\in \bar{\Omega}}\lim_{\stackrel{x\rightarrow y}{x\in \Omega}}|x-y|^{\mu+p}V_2(x)<\infty$.
    \item[$(H_2)$] $1<q<p, \Omega$ is unbounded  (in particular, $\Omega=\R^N$),
     $$ \sup_{y\in \overline{\Omega\cap B_R(0)}}\lim_{\stackrel{x\rightarrow y}{x\in \Omega}}|x-y|^{\mu+p}V_2(x)<\infty\;\hbox{for any fixed}\; R>0$$ and
        $$ \lim_{R\rightarrow \infty}\int_{\{x\in \Omega:|x|>R\}} \big(V_2(x)\big)^{\frac{p}{p-q}}|x|^{\frac{(\mu+p)q}{p-q}}dx\rightarrow 0.$$
\item[$(H_3)$]$p\leq q. $   For any $\displaystyle y\in \bar{\Omega}, \lim_{\stackrel{x\rightarrow y}{x\in \Omega}}|x-y|^{\bar{\sigma}}V_2(x)=0$ and $\displaystyle \lim_{\stackrel{|x|\rightarrow \infty}{x\in\Omega}}|x|^{\bar{\sigma}}V_2(x)=0$, where
    $$\bar{\sigma}:=N-\frac{q}{p}(N-\mu-p)\in(\frac{N\mu}{N-p}, \mu+p].$$
    \end{itemize}
\end{itemize}
 Here comes our another main theorem:\\

\bt\lab{2014-9-7-th3}  Assume ${\bf (H)}$.
 \begin{itemize}
 \item[(1)] If $q=p$,  then  the equation (\ref{2014-9-7-xe2}) has a sequence of  eigenfunctions   $\{\varphi_n\}$,  the corresponding  eigenvalues $\{ \lambda_n\}$  satisfying  $\lambda_n\rightarrow \infty$ as $n\rightarrow \infty$.
 \item[(2)] If  $q>p$,   then  for any positive fixed $\lambda$, \eqref{2014-9-7-xe2} possesses  a sequence of solutions $\{v_n\}$ such that $0<c_1\leq c_2\leq \cdots\leq c_n\rightarrow \infty$ as $n\rightarrow \infty$.
 \item[(3)] If  $q<p$,   then  for any positive fixed $\lambda$, \eqref{2014-9-7-xe2} possesses  a sequence of solutions $\{v_n\}$ such that $-\infty<c_1\leq c_2\leq\cdots\leq c_n\rightarrow 0$ as $n\rightarrow \infty$.
  \end{itemize}
 Where
 $$c_n:=\Phi(v_n):=\frac{1}{p}\int_\Omega \frac{|\nabla v_n|^p}{|x|^\mu}dx-\frac{1}{q}\lambda\int_\Omega V(x)|v_n|^qdx.$$
\et

\br  When $\mu=0$ and $q=p$,  Theorem \ref{2014-9-7-th3}  was established by Szulkin-Willem   \cite{SzulkinWillem.1999}.
\er

\s{Preliminaries}
\renewcommand{\theequation}{2.\arabic{equation}}
\renewcommand{\theremark}{2.\arabic{remark}}
\renewcommand{\thedefinition}{2.\arabic{definition}}
Firstly, based  on the  original Rellich-Kondrachov theorem and the general weighted Hardy inequality (\ref{2013-9-05-e8}), we obtain the following result through a   transformation.
\bl\lab{2014-9-6-l1}
Assume that $1\leq p<N, N-p-\mu>0$ and  that $\{u_n\}\subset W_{\mu}^{1,p}(\Omega)$ is bounded. Then there exists some $u\in W_{\mu}^{1,p}(\Omega)$ and,  up to a subsequence,  $u_n\rightarrow u$ a.e.  in $\Omega$.
\el
\bp
For the case of $\mu=0$, it  can be easily obtained by the original Rellich-Kondrachov theorem and the diagonal trick. Next, we only consdier that $\mu\neq 0$.
Denote $\bar{\mu}:=\frac{\mu}{p}$, note that
$$\nabla \big(\frac{u_n}{|x|^{\bar{\mu}}}\big)=\frac{\nabla u_n}{|x|^{\bar{\mu}}}-\bar{\mu}|x|^{-\bar{\mu}-2}xu_n,$$
it follows that
$$\big|\nabla (\frac{u_n}{|x|^{\bar{\mu}}})\big|\leq \frac{|\nabla u_n|}{|x|^{\bar{\mu}}}+|\bar{\mu}|\frac{|u_n|}{|x|^{\bar{\mu}+1}}.$$
Recalling that for $p>0$, $
|x+y|^p\leq \max\{1, 2^{p-1}\}(|x|^p+|y|^p). $
Hence, combining with the general weighted Hardy inequality (\ref{2013-9-05-e8}) due to the fact of $N-p-\mu>0$, we have
\begin{align}\lab{2014-9-6-xe2}
&\int_\Omega \big|\nabla (\frac{u_n}{|x|^{\bar{\mu}}})\big|^pdx\nonumber\\
& \leq 2^{p-1}\Big[\int_\Omega \big|\frac{|\nabla u_n|}{|x|^{\bar{\mu}}}\big|^pdx+|\bar{\mu}|^p\int_\Omega \big|\frac{|u_n|}{|x|^{\bar{\mu}+1}}\big|^pdx\Big]\nonumber\\
& \leq  C(p,\mu) \int_\Omega \big|\frac{|\nabla u_n|}{|x|^{\bar{\mu}}}\big|^pdx=C(p,\mu) \int_\Omega \frac{|\nabla u_n|^p}{|x|^\mu}dx.
\end{align}
Thus, $\displaystyle\{\frac{u_n}{|x|^{\bar{\mu}}}\}$ is a bounded sequence in $W^{1,p}(\Omega)$. It follows  from  the well-known Rellich-Kondrachov compactness theorem and the standard diagonal trick, we obtain that, up to a subsequence, $\displaystyle\frac{u_n}{|x|^{\bar{\mu}}}\rightarrow \frac{u}{|x|^{\bar{\mu}}}$ a.e. in $\Omega$.
Then it is natural to see that $u_n\rightarrow u$ a.e.  in $\Omega$.
It follows from the Fatou's Lemma that $u\in W_{\mu}^{1,p}(\Omega)$.
\ep

Now, we can prove the weighted Rellich-Kondrachov compactness theorem:

\bt\lab{2014-9-7-th1}
Assume $\Omega\subset \R^N$ is a bounded open Lipschitz domain. Suppose that $1\leq p<N, -\infty<\mu<N-p$, then the embedding
$$W_{\mu}^{1,p}(\Omega)\hookrightarrow L^q(\Omega, \frac{dx}{|x|^s})$$ is compact if
$$\begin{cases}
\frac{N\mu}{N-p}\leq s\leq \mu+p,\\
1\leq q<p^*(s, \mu),
\end{cases}\;\hbox{or}\;\begin{cases}
s<\frac{N\mu}{N-p},\\
1\leq q<p^*.
\end{cases}$$
Moreover, if $s\geq \max\{0,\frac{N\mu}{N-p}\}$,  the  conclusion  is still valid when  domain $\Omega$ is unbounded but with finite Lebesgue measure.  \et

\bp  We firstly consider the case  that $\Omega$ is bounded.  Assume that $\displaystyle \sup_{n}\|u_n\|<\infty$.
By Lemma \ref{2014-9-6-l1}, without loss of generality, we may assume that $u_n\rightarrow u$ a.e. in $\Omega$ for some $u\in W_{\mu}^{1,p}(\Omega)$.
Since $s<N$ and $\Omega$ is bounded, it is easy to see that $\nu(\Omega)<\infty$, where the new measure $d\nu:=\frac{dx}{|x|^s}$ and $\nu\big|_\Omega$ is absolutely continuous  with respect to   the usual Lebesgue measure $L$.
If $\frac{N\mu}{N-p}\leq s\leq \mu+p$, by the Hardy-Sobolev inequality (\ref{2013-9-05-e5}), we also have that
\be\lab{2014-9-7-xe3}
\sup_{n}\int_\Omega |u_n|^{p^*(s,\mu)}d\nu<\infty.
\ee
Then by the H\"{o}lder inequality related to the measure $\nu$, for any subset $\Lambda\subset \Omega$, since $1\leq q<p^*(s,\mu)$, we have
\begin{align}\lab{2014-9-6-xe3}
&\int_{\Lambda}|u_n-u|^qd\nu\nonumber\\
&\leq \Big(\int_{\Lambda}|u_n-u|^{p^*(s,\mu)}d\nu\Big)^{\frac{q}{p^*(s,\mu)}}
\Big(\nu(\Lambda)\Big)^{\frac{p^*(s, \mu)-q}{p^*(s,\mu)}}\nonumber\\
&\leq  C\Big(\nu(\Lambda)\Big)^{\frac{p^*(s, \mu)-q}{p^*(s,\mu)}}.
\end{align}
Recalling that $\nu\big|_\Omega$ is absolutely continuous, we have $\nu(\Lambda)\rightarrow 0$ as $L(\Lambda)\rightarrow 0$. Hence,
\be\lab{2014-9-23-e1}
\int_{\Lambda}\frac{|u_n-u|^q}{|x|^s}dx\rightarrow 0 \;\hbox{as the measure}\;L(\Lambda)\rightarrow 0\;\hbox{uniformly for all $n$}.
\ee
Since $\Omega$ is bounded and $u_n\rightarrow u$ a.e. in $\Omega$, applying the Egoroff Theorem and the above conclusion (\ref{2014-9-23-e1}), we see that,  up to a subsequence,
$$\int_\Omega \frac{|u_n-u|^q}{|x|^s}dx\rightarrow 0\;\hbox{as}\;n\rightarrow \infty.$$
When  $s<\frac{N\mu}{N-p}, 1\leq q<p^*,$   we  denote $s_0:=\frac{N\mu}{N-p}$. Note  that $p^*(s_0,\mu)=\frac{Np}{N-p}=p^*$. By the above arguments, if $1\leq q<p^*$, we have
\be\lab{2014-9-30-xe1}
u_n\rightarrow u\;\hbox{strongly in}\;L^q(\Omega, \frac{dx}{|x|^{s_0}}).
\ee
Hence,
\begin{align}\lab{2014-9-30-xe2}
\int_\Omega \frac{|u_n-u|^q}{|x|^s}dx=\int_\Omega \frac{|u_n-u|^q}{|x|^{s_0}} |x|^{s_0-s}dx
\leq C(\Omega)\int_\Omega \frac{|u_n-u|^q}{|x|^{s_0}}dx
\rightarrow 0 \end{align}
as $ n\rightarrow \infty.$ \ep

To consider the case with unbounded $\Omega$, we  insert   the    definition of tightness which can be found in \cite{Lions.1985a,Lions.1985}.

\bd\lab{2013-8-4-def1}
Assume $\{\rho_k\}$ is a bounded sequence in $L^1(\R^N)$ and $\rho_k\geq 0$ satisfies
\be\lab{2013-8-4-e1}
\|\rho_k\|_{L^1}=\lambda+o(1), \;\lambda>0.
\ee
Then we call this sequence $\{\rho_k\}$ is a tight sequence if $\forall\;\varepsilon>0, \exists\;R>0$ such that
\be\lab{2013-8-4-e2}
\int_{|x|\geq R}\rho_k(x)dx<\varepsilon,\;\forall\;k\geq 1.
\ee
We call $u_k$ is a $L^p$ tight sequence, if $|u_k|^p$ is a tight sequence.  For the  convenience, the definition is still valid in the   current  paper when $\lambda=0$ in (\ref{2013-8-4-e1}).
\ed

\vskip 0.03in
\noindent
{\bf Completion of the proof of Theorem \ref{2014-9-7-th1}.}   We assume that $s\geq \max\{0, \frac{N\mu}{N-p}\}$ and $L(\Omega)<\infty$. Note  that in this case, $p^*(s,\mu)\leq p^*$.  Based  on the results above, we only need to show that $\displaystyle\{\frac{|u_n|^q}{|x|^s}\}$ is a tight sequence. It is still satisfying  that $\nu(\Omega)<\infty$ and $\nu\big|_\Omega$ is absolutely continuous, where $\nu=\frac{dx}{|x|^s}$.  Hence,
$$L(\Omega \cap B_R^c(0))\rightarrow 0\;\hbox{as}\;R\rightarrow \infty.$$
and then
\be\lab{2014-9-6-xe4}
\nu(\Omega \cap B_R^c(0))\rightarrow 0\;\hbox{as}\;R\rightarrow \infty.
\ee
Apply the H\"{o}lder inequality on the domain $\Omega \cap B_R^c(0)$, it follows  from (\ref{2014-9-7-xe3}) and (\ref{2014-9-6-xe4}) that
$\displaystyle\{\frac{|u_n|^q}{|x|^s}\}$ is a tight sequence.\hfill$\Box$

\newpage

\s{The existence of extremal functions  for a family of CKN  inequalities}
\renewcommand{\theequation}{3.\arabic{equation}}
\renewcommand{\theremark}{3.\arabic{remark}}
\renewcommand{\thedefinition}{3.\arabic{definition}}

Firstly, consider the following problem:
\be\lab{2014-9-9-xe2}
-div\big(\frac{|\nabla u|^{p-2}\nabla u}{|x|^\mu}\big)+\lambda \frac{|u|^{q-2}u}{|x|^\sigma}=\frac{|u|^{r-2}u}{|x|^s}\;\hbox{ in }\;\Omega, \; \lambda>  0, \; u\in W_{\mu}^{1,p}(\Omega),\\
\ee
where $p>1, \max\{\sigma, s\}<\mu+p<N, q\geq 1, r\geq 1$ and $\Omega$ is an open Lipschitz domain of $\R^N$.  We introduce the Sobolev space \be\lab{2014-9-9-xe3-0} E:=W_{\mu}^{1,p}(\Omega)\cap L^q(\Omega, \frac{dx}{|x|^\sigma})\ee
which is equipped with the norm
\be\lab{2014-9-9-xe3}
\|u\|_E:=\|u\|+\lambda |u|_{q,\sigma},
\ee
where $\|u\|$ is the   norm in $W_{\mu}^{1,p}(\Omega) $ defined by (\ref{zz=z}).
Then we see that $E$ is a closed subspace of $W_{\mu}^{1,p}(\Omega)$.

\bl\lab{2014-9-19-l1}
Assume that $\displaystyle p>1, \max\{\sigma, s\}< \mu+p<N$, $\displaystyle  1\leq q<\min\{p^*,p^*(\sigma, \mu)\}$, $\displaystyle \max\{\frac{p(\sigma-s)}{N-\mu-p}+q, \frac{\sigma-s}{N-\sigma}q+q\}<r<\min\{p^*, p^*(s, \mu)\}$ and  that
$$
\begin{cases}
p(s-\sigma)+q(\mu+p-s)\leq r(\mu+p-\sigma),\\
(Np-Nr+pr)(s-\sigma)\geq (N\mu-Ns+ps)(r-q),
\end{cases}
$$
 then there exists some constant $C(p,q,r,\mu,\sigma,s)>0$ such that the CKN   inequality \eqref{2014-6-5-xe4} holds true,  i.e.,
\be\lab{2014-9-19-e1}
\big(\int_\Omega \frac{|u|^r}{|x|^s}dx\big)^{\frac{1}{r}}\leq C(p,q,r,\mu,\sigma,s)\big(\int_{\Omega}\frac{|\nabla u|^p}{|x|^\mu}dx\big)^{\frac{a}{p}}\big(\int_\Omega \frac{|u|^q}{|x|^\sigma}\big)^{\frac{1-a}{q}}
\ee
for all $u\in E:=W_{\mu}^{1,p}(\Omega)\cap L^q(\Omega, \frac{dx}{|x|^\sigma})$,
where $$a=\frac{[(N-\sigma)r-(N-s)q]p}{[(N-\sigma)p-(N-\mu-p)q]r}\in (0,1).$$
Moreover, if $r\geq 1$, $E\hookrightarrow L^r(\Omega,\frac{dx}{|x|^s})$ is a continuous embedding.
\el

\bp
Let
\be\lab{2014-9-19-xe1}
\begin{cases}
r_1:=\frac{\big[p(N-s)-r(N-\mu-p)\big]q}{p(N-\sigma)-q(N-\mu-p)},\\
r_2:=\frac{\big[r(N-\sigma)-q(N-s)\big]p}{p(N-\sigma)-q(N-\mu-p)},\\
s_1:=\frac{\big[p(N-s)-r(N-\mu-p)\big]\sigma}{p(N-\sigma)-q(N-\mu-p)},\\
s_2:=\frac{Np(s-\sigma)+(r\sigma-qs)(N-\mu-p)}{p(N-\sigma)-q(N-\mu-p)},\\
\bar{\sigma}:=\frac{Np(s-\sigma)+(r\sigma-qs)(N-\mu-p)}{p(s-\sigma)+(r-q)(N-\mu-p)},
\end{cases}
\ee
then a direct calculation shows that
\be\lab{2014-9-19-xe2}
\begin{cases}
0<r_1<q, r_2>0, r_1+r_2=r,\\
s_1+s_2=s,\frac{N\mu}{N-p}\leq \bar{\sigma}\leq \mu+p<N,\\
\frac{q s_1}{r_1}=\sigma,\\
\frac{q s_2}{q-r_1}=\bar{\sigma},\\
\frac{q r_2}{q-r_1}=p^*(\bar{\sigma}, \mu):=\frac{p(N-\bar{\sigma})}{N-\mu-p}\in [p, p^*].
\end{cases}
\ee
Thus, by the H\"{o}lder inequality and Hardy-Sobolev inequality, we have
\begin{align}\lab{2014-9-19-xe3}
&\int_{\Omega}\frac{|u|^r}{|x|^s}dx=\int_{\Omega}\frac{|u|^{r_1}}{|x|^{s_1}}\cdot \frac{|u|^{r_2}}{|x|^{s_2}}dx\nonumber\\
& \leq \Big(\int_{\Omega}\frac{|u|^q}{|x|^\sigma}dx\Big)^{\frac{r_1}{q}}
\Big(\int_{\Omega}\frac{|u|^{p^*(\bar{\sigma}, \mu)}}{|x|^{\bar{\sigma}}}dx\Big)^{\frac{q-r_1}{q}}\nonumber\\
& \leq  C(\bar{\sigma})\Big(\int_{\Omega}\frac{|u|^q}{|x|^\sigma}dx\Big)^{\frac{r_1}{q}}
\Big(\int_{\Omega}\frac{|\nabla u|^p}{|x|^{\mu}}dx\Big)^{\frac{p^*(\bar{\sigma},\mu)}{p}\frac{q-r_1}{q}}.
\end{align}
Note that $0<1-a=\frac{r_1}{r}<1$ since $0<r_1<r$, we also have
$\displaystyle \frac{p^*(\bar{\sigma},\mu)(q-r_1)}{qr}=a.$
Hence, we obtain that there exists some $C(p,q,r,\mu,\sigma,s)>0$ such that
(\ref{2014-9-19-e1}) is satisfied
for all $u\in W_{\mu}^{1,p}(\Omega)\cap L^q(\Omega, \frac{dx}{|x|^\sigma})$.
Finally, if $r\geq 1$, by the Young inequality, we have
\be\lab{2014-9-19-xe4}
|u|_{r,s}\leq \max\big\{\frac{1}{(1-a)\lambda}, \frac{1}{a}\big\}C(p,q,r,\mu,\sigma,s) \|u\|_E,
\ee
where
$\|u\|_E$ is defined by \eqref{2014-9-9-xe3}.
Thus, $E\hookrightarrow L^r(\Omega,\frac{dx}{|x|^s})$ is a continuous embedding.
\ep

\bl\lab{2014-9-19-l2}
Under the assumptions of Lemma \ref{2014-9-19-l1}, if furthermore $s>0$ and the condition \eqref{2014-9-30-bu2} strictly holds, i.e.,
 \be\lab{2014-9-30-bu3}
\begin{cases}
p(s-\sigma)+q(\mu+p-s)< r(\mu+p-\sigma)\\
(Np-Nr+pr)(s-\sigma)> (N\mu-Ns+ps)(r-q)
\end{cases},
\ee
 then any bounded sequence $\{u_n\}$ of $E:=W_{\mu}^{1,p}(\Omega)\cap L^q(\Omega, \frac{dx}{|x|^\sigma})$ satisfying  that $\{\frac{|u_n|^r}{|x|^s}\}$ is a tight sequence. In particular, if $r\geq 1$, the embedding
$$E\hookrightarrow L^r(\Omega,\frac{dx}{|x|^s})$$
is compact.
\el
\bp
Let $\{u_n\}\subset E$ be a bounded sequence. Since $s>0$, by the continuity,  we can take some $0<\bar{s}<s$ close to $s$ such that the assumptions of Lemma \ref{2014-9-19-l1} still hold after replacing $s$ by $\bar{s}$. Thus
$\sup_{n} |u_n|_{r,\bar{s}}\leq C.$
Then it follows that
\begin{align}\lab{2014-9-19-xe5}
&\int_{|x|>R}\frac{|u_n|^r}{|x|^s}dx=\int_{|x|>R}\frac{1}{|x|^{s-\bar{s}}}\frac{|u_n|^r}{|x|^{\bar{s}}}dx\nonumber\\
&< R^{\bar{s}-s}|u_n|_{r,\bar{s}}^{r}\nonumber\\
& \rightarrow 0\;\hbox{uniformly for all $n$ as}\;R\rightarrow \infty.
\end{align}
Hence, $\{\frac{|u_n|^r}{|x|^s}\}$ is a tight sequence.
Recalling the  Theorem \ref{2014-9-7-th1} for  the case bounded domain, it is easy to prove that the embedding
$E\hookrightarrow L^r(\R^N,\frac{dx}{|x|^s})$ is compact.
\ep

\bc\lab{2014-10-1-cro1}
Let $\sigma=0, 1<p=q<N$ and one of the following holds:
\begin{itemize}
\item[(i)]$\mu=0, 0<s<p<N, p\leq r<\frac{p(N-s)}{N-p}$;
\item[(ii)]$\mu>0, \frac{N\mu(r-p)}{p^2}< s<\mu+p<N, p\leq r<\min\{\frac{pN}{N-p}, \frac{p(N-s)}{N-\mu-p}\}$;
\item[(iii)] $\mu<0, 0<s<\mu+p<N, p\leq r< \frac{p(N-s)}{N-\mu-p}$,
\end{itemize}
then $$W_{\mu}^{1,p}(\Omega)\cap L^p(\Omega)\hookrightarrow L^r(\Omega, \frac{dx}{|x|^s})$$
is a compact embedding.
\ec
\bp
We note that for case $r=p$, we can apply the similar arguments as the proof of Lemma \ref{2014-9-19-l2}. And other cases are straight-forward results of Lemma \ref{2014-9-19-l2}.
\ep

\br\lab{2014-9-19-r1}
It is well known that $H^{1,p}(\R^N)\hookrightarrow L^r(\R^N)$ is a continuous embedding for $r\in [p, p^*)$ but not compact. Take $\mu=0,1<p<N,0<s<p$, then by (i) of Corollary \ref{2014-10-1-cro1} we see that $H^{1,p}(\R^N)\hookrightarrow L^r(\R^N, \frac{dx}{|x|^s})$ is a  compact embedding for $r\in [p, p^*(s))$.
\er

\bl\lab{2014-9-21-l1}
Assume that $a>0,b>0,A>0,B>0$ are fixed, let $g(t):=t^aA+t^{-b}B$, then
$$\inf_{t>0}g(t)=g(t_0)=\frac{a+b}{a} (\frac{b}{a})^{\frac{-b}{a+b}}A^{\frac{b}{a+b}}B^{\frac{a}{a+b}},$$
where
$t_0=(\frac{bB}{aA})^{\frac{1}{a+b}}.$
\el
\bp
It is  a direct computation. \ep

\bt\lab{2014-9-21-th1}
Let $\Omega$ be a cone, that is,  $\Omega=\frac{\Omega}{t}$ for any $t>0$ (in particular,  $\Omega=\R^N$). Assume that $\displaystyle p>1, \max\{\sigma, s\}< \mu+p<N$, $\displaystyle  1\leq q<\min\{p^*,p^*(\sigma, \mu)\}$, $\displaystyle \max\{\frac{p(\sigma-s)}{N-\mu-p}+q, \frac{\sigma-s}{N-\sigma}q+q\}<r<\min\{p^*, p^*(s, \mu)\}$ and
$$
\begin{cases}
p(s-\sigma)+q(\mu+p-s)\leq r(\mu+p-\sigma)\\
(Np-Nr+pr)(s-\sigma)\geq (N\mu-Ns+ps)(r-q)
\end{cases},
$$
then
\begin{align}\lab{2014-9-21-xe1}
&\inf_{u\in M}I(u)\nonumber\\
=&\frac{(\mu+p-\sigma)r+(p-q)(N-s)}{p(N-s)-(N-\mu-p)r}\big[\frac{(N-\sigma)r-(N-s)q}{p(N-s)-(N-\mu-p)r}\big]^{\frac{(N-s)q-(N-\sigma)r}{(p-q)(N-s)+(\mu+p-\sigma)r}}\nonumber\\
&\lambda^{\frac{(N-s)p-(N-\mu-p)r}{(\mu+p-\sigma)r+(p-q)(N-s)}}\inf_{u\in M}\Big\{\|u\|^{\frac{p[(N-\sigma)r-(N-s)q]}{(\mu+p-\sigma)r+(p-q)(N-s)}}|u|_{q,\sigma}^{\frac{q[(N-s)p-(N-\mu-p)r]}{(\mu+p-\sigma)r+(p-q)(N-s)}}\Big\},
\end{align}
where
\be\lab{2014-9-19-we1}
M:=\{u\in E:     \int_\Omega \frac{|u|^r}{|x|^s}dx=1\},
\ee
and
\be\lab{2014-9-26-xe1}
I(u):=\int_\Omega\frac{|\nabla u|^p}{|x|^\mu}dx+\lambda \int_\Omega\frac{|u|^q}{|x|^\sigma}dx.
\ee
\et
\bp
For $t>0$, define a mapping $T_t(u)=u_t:=t^{\frac{N-s}{r}}u(tx)$, then it is easy to check that
$$\int_\Omega \frac{|u|^r}{|x|^s}dx\equiv \int_\Omega \frac{|u_t|^r}{|x|^s}dx, t>0.$$
Hence, $M$ is invariant under the transformation $T_t$.
We note that $$\frac{p(N-s)-(N-\mu-p)r}{r}>0, \frac{(N-\sigma)r-(N-s)q}{r}>0.$$ Hence, by Lemma \ref{2014-9-21-l1}, we have
\begin{align*}
&\inf_{t>0}I(u_t)\\
=&\inf_{t>0}\Big(t^{\frac{p(N-s)-(N-\mu-p)r}{r}}\int_\Omega \frac{|\nabla u(x)|^p}{|x|^\mu}dx+t^{-\frac{r(N-\sigma)-q(N-s)}{r}}\lambda\int_\Omega \frac{|u(x)|^q}{|x|^\sigma}dx\Big)\\
=&\frac{(\mu+p-\sigma)r+(p-q)(N-s)}{p(N-s)-(N-\mu-p)r}\big[\frac{(N-\sigma)r-(N-s)q}{p(N-s)-(N-\mu-p)r}\big]^{\frac{(N-s)q-(N-\sigma)r}{(p-q)(N-s)+(\mu+p-\sigma)r}}\nonumber\\
&\lambda^{\frac{(N-s)p-(N-\mu-p)r}{(\mu+p-\sigma)r+(p-q)(N-s)}}\|u\|^{\frac{p[(N-\sigma)r-(N-s)q]}{(\mu+p-\sigma)r+(p-q)(N-s)}}|u|_{q,\sigma}^{\frac{q[(N-s)p-(N-\mu-p)r]}{(\mu+p-\sigma)r+(p-q)(N-s)}}\nonumber\\
=:&C^*(p,q,r,\mu,\sigma,s,\lambda)\|u\|^{\frac{p[(N-\sigma)r-(N-s)q]}{(\mu+p-\sigma)r+(p-q)(N-s)}}|u|_{q,\sigma}^{\frac{q[(N-s)p-(N-\mu-p)r]}{(\mu+p-\sigma)r+(p-q)(N-s)}}.
\end{align*}
Thereby, we prove this Lemma.
\ep

\br\lab{2014-9-21-r1}
Define
\begin{align}\lab{2014-9-21-xe2}
\lambda^*(p,q,r,\mu,\sigma,s) := &\Big\{\frac{p(N-s)-(N-\mu-p)r}{(\mu+p-\sigma)r+(p-q)(N-s)}\Big\}^{\frac{(\mu+p-\sigma)r+(p-q)(N-s)}{(N-s)p-(N-\mu-p)r}}\nonumber\\
&\cdot\Big\{\frac{(N-\sigma)r-(N-s)q}{p(N-s)-(N-\mu-p)r}\Big\}^{\frac{(N-\sigma)r-(N-s)q}{p(N-s)-(N-\mu-p)r}},
\end{align}
then we have
$$C^*\big(p,q,r,\mu,\sigma,s,\lambda^*(p,q,r,\mu,\sigma,s)\big)\equiv 1.$$
For  the simplicity,  if there exists no misunderstanding,  we will write $$\lambda^*=\lambda^*(p,q,r,\mu,\sigma,s);\quad I^*(u)=\int_\Omega \frac{|\nabla u|^p}{|x|^\mu}dx+\lambda^* \int_\Omega\frac{|u|^q}{|x|^\sigma}dx. $$ \er

\bc\lab{2014-9-21-cro1}
Let $\Omega$ be a cone (i.e., $\Omega=\frac{\Omega}{t}$ for any $t>0$.  In particular,  $\Omega=\R^N$).  Assume that $\displaystyle p>1, \max\{\sigma, s\}< \mu+p<N$, $\displaystyle  1\leq q<\min\{p^*,p^*(\sigma, \mu)\}$, $\displaystyle \max\{\frac{p(\sigma-s)}{N-\mu-p}+q, \frac{\sigma-s}{N-\sigma}q+q\}<r<\min\{p^*, p^*(s, \mu)\}$ and
$$
\begin{cases}
p(s-\sigma)+q(\mu+p-s)\leq r(\mu+p-\sigma)\\
(Np-Nr+pr)(s-\sigma)\geq (N\mu-Ns+ps)(r-q)
\end{cases},
$$
then the sharp constant of (\ref{2014-9-19-e1})
$$C(p,q,r,\mu,\sigma,s)=\big(\frac{1}{\rho}\big)^{\frac{(\mu+p-\sigma)q+(p-r)(N-s)}{q[(N-\sigma)p-(N-\mu-p)r]}},$$
where
$\rho:=\inf_{u\in M}I^*(u).$
\ec
\bp
For any $u\in M$, by Theorem \ref{2014-9-21-th1}, we have
\be\lab{2014-9-21-xe3}
\inf_{t>0}I^*(u_t)=\|u\|^{\frac{p[(N-\sigma)r-(N-s)q]}{(\mu+p-\sigma)r+(p-q)(N-s)}}|u|_{q,\sigma}^{\frac{q[(N-s)p-(N-\mu-p)r]}{(\mu+p-\sigma)r+(p-q)(N-s)}}.
\ee
It follows that
\be\lab{2014-9-21-xe4}
\rho:=\inf_{u\in M}I^*(u)=\inf_{u\in M}\|u\|^{\frac{p[(N-\sigma)r-(N-s)q]}{(\mu+p-\sigma)r+(p-q)(N-s)}}|u|_{q,\sigma}^{\frac{q[(N-s)p-(N-\mu-p)r]}{(\mu+p-\sigma)r+(p-q)(N-s)}}.
\ee
Note  that
\begin{align*}
&\frac{p[(N-\sigma)r-(N-s)q]}{(\mu+p-\sigma)r+(p-q)(N-s)}+\frac{q[(N-s)p-(N-\mu-p)r]}{(\mu+p-\sigma)r+(p-q)(N-s)}\\
& = \frac{r[(N-\sigma)p-(N-\mu-p)q]}{(\mu+p-\sigma)r+(p-q)(N-s)},
\end{align*}
then we have that
$$\rho|u|_{r,s}^{\frac{r[(N-\sigma)p-(N-\mu-p)q]}{(\mu+p-\sigma)r+(p-q)(N-s)}}\leq \|u\|^{\frac{p[(N-\sigma)r-(N-s)q]}{(\mu+p-\sigma)r+(p-q)(N-s)}}|u|_{q,\sigma}^{\frac{q[(N-s)p-(N-\mu-p)r]}{(\mu+p-\sigma)r+(p-q)(N-s)}} $$
for all $u\in E$,
it follows that
\be\lab{2014-9-21-xe5}
|u|_{r,s}\leq \big(\frac{1}{\rho}\big)^{\frac{(\mu+p-\sigma)r+(p-q)(N-s)}{r[(N-\sigma)p-(N-\mu-p)q]}}\|u\|^{a}|u|_{q,\sigma}^{1-a}\;\hbox{for all}\;u\in E.
\ee
Note that the above processes   are  reversible,  the Corollary is proved.
\ep

\br\lab{2014-9-21-r2}
Under the assumptions of Corollary \ref{2014-9-21-cro1}, the sharp constant of inequality (\ref{2014-9-19-e1}) can be achieved if and only if $\rho$ can be reached.
\er

\vskip 0.1in
Next, let us assume $r\geq 1$ and consider the following  minimizing problem.

\bl\lab{2014-9-19-l3}
Assume that $\displaystyle p>1, \max\{\sigma, s\}< \mu+p<N$, $\displaystyle  r\geq 1, 1\leq q<\min\{p^*,p^*(\sigma, \mu)\}$, $\displaystyle \max\{\frac{p(\sigma-s)}{N-\mu-p}+q, \frac{\sigma-s}{N-\sigma}q+q\}<r<\min\{p^*, p^*(s, \mu)\}, s>0$ and
$$
\begin{cases}
p(s-\sigma)+q(\mu+p-s)< r(\mu+p-\sigma)\\
(Np-Nr+pr)(s-\sigma)> (N\mu-Ns+ps)(r-q)
\end{cases},
$$
then   $\rho$  can be achieved by some minimizer $u\in M$.
Furthermore,  if $q>1, r>1$, the minimizer is a ground state solution to  the following problem:
\begin{align}\lab{2014-9-21-wwe1}
-div(\frac{|\nabla u|^{p-2}\nabla u}{|x|^\mu})+\frac{q\lambda^*}{p}\frac{|u|^{q-2}u}{|x|^\sigma}
= \frac{r[p(N-\sigma)-(N-\mu-p)q] \rho}{p[(\mu+p-\sigma)r+(p-q)(N-s)]}  \frac{|u|^{r-2}u}{|x|^s},  \end{align}
for $u\in E,$
where $\rho$ is defined in Corollary \ref{2014-9-21-cro1}.
\el
\bp
Obviously,  $\lambda^*>0$.  Let $\{u_n\}$ be a minimizing sequence of $\rho$ in $M$, i.e., $|u_n|_{r,s}\equiv 1$ and $I^*(u_n)\rightarrow \rho$.
By (\ref{2014-9-19-xe4}),   $\rho>0$. Further,  $\{u_n\}$ is bounded in $E$. By Lemma \ref{2014-9-19-l2}, up to a subsequence if necessary, we may assume that $u_n\rightarrow u$ in $L^r(\Omega,\frac{dx}{|x|^s})$. Hence, $u\in M.$
We also note that $\sup_{n\geq 1}\|u_n\|<\infty$, then by Lemma \ref{2014-9-6-l1}, we may assume that $u_n\rightarrow u$ a.e.  in $\Omega$ since $p>1$. Follows  from  the Fatou's Lemma, we have $
I^*(u)\leq \liminf_{n\rightarrow \infty}I^*(u_n)=\rho.
$
On the other hand, by the definition of $\rho$, we have
$
I^*(u)\geq \rho\;\hbox{since}\;u\in M.
$
Hence, $u$ is a minimizer.   Let $u$ be an  extremal function.  If $q>1, r>1$, then there exists some Lagrange multiplier $\tilde{\lambda}$ such that
\be\lab{2014-9-21-wwe2}
-p\;div(\frac{|\nabla u|^{p-2}\nabla u}{|x|^\mu})+q\lambda^*\frac{|u|^{q-2}u}{|x|^\sigma}=\tilde{\lambda}\frac{|u|^{r-2}u}{|x|^s}.
\ee
Testing by $u$, we obtain that
\be\lab{2014-9-21-wwe5}
p\|u\|^p+q\lambda^*|u|_{q,\sigma}^{q}=\tilde{\lambda}.
\ee
Recalling that
\be\lab{2014-9-21-wwe6}
\|u\|^p+\lambda^*|u|_{q,\sigma}^{q}=\rho
\ee
and by Lemma \ref{2014-9-21-l1}, we see that
\be\lab{2014-9-21-wwe7}
\frac{p(N-s)-(N-\mu-p)r}{r}\|u\|^p=\frac{r(N-\sigma)-q(N-s)}{r}\lambda^*|u|_{q,\sigma}^{q}.
\ee
Combine  (\ref{2014-9-21-wwe5}), (\ref{2014-9-21-wwe6}) and (\ref{2014-9-21-wwe7}), we obtain that
\be\lab{2014-9-21-wwe8}
\tilde{\lambda}=\frac{r[p(N-\sigma)-(N-\mu-p)q]}{(\mu+p-\sigma)r+(p-q)(N-s)}\rho.
\ee
Hence, the minimizer is a ground state solution to  the   equation  (\ref{2014-9-21-wwe1}).
 \ep

\noindent{\bf Proof of Theorem \ref{2014-10-14-mainth1}. }  It is a straightforward consequence
of Corollary \ref{2014-9-21-cro1}  and Lemma \ref{2014-9-19-l3}.  \hfill  $\Box$

\s{An application}
\renewcommand{\theequation}{4.\arabic{equation}}
\renewcommand{\theremark}{4.\arabic{remark}}
\renewcommand{\thedefinition}{4.\arabic{definition}}
In this section, we will study the problem (\ref{2014-9-7-xe2}) as an application of  the previous theorem.
Based on the Hardy-Sobolev inequality (\ref{2013-9-05-e5}),  firstly, we get the following
weakly continuous functional, which   largely generalizes  the corresponding result in  \cite[Lemma 2.13]{Willem.1996}.

\bl\lab{2014-9-8-l1}
If $1<p<N, \mu+p<N, q<p^*:=\frac{pN}{N-p}$ and there exists some $\frac{N\mu}{N-p}\leq\eta< \min\{\mu+p, N-\frac{q}{p}(N-\mu-p)\}$ such that $\big|a(x)\big|^{\frac{p^*(\eta, \mu)}{p^*(\eta,\mu)-q}}|x|^{\frac{q\eta}{p^*(\eta,\mu)-q}}\in L^1(\Omega)$, then the functional $\chi:W_{\mu}^{1,p}(\Omega)\rightarrow \R $ defined by
$$\chi(u)= \int_\Omega a(x)|u|^qdx$$
is weakly continuous. In particular,  when $\mu\leq 0<\mu+p<N, q=p$ and $ a(x)\in L^{\frac{N}{\mu+p}}(\Omega)$, the result holds.
\el
\bp
By $\frac{N\mu}{N-p}\leq\eta< \min\{\mu+p, N-\frac{q}{p}(N-\mu-p)\}$,  we have $p^*\geq p^*(\eta,\mu)>\max\{p,q\}$.
Then by  the  H\"{o}lder inequality, for any $u\in W_{\mu}^{1,p}(\Omega)$, we have
\begin{align*}
&\int_\Omega |a(x)||u|^qdx=\int_\Omega \big|a(x)|x|^{\frac{q\eta}{p^*(\eta,\mu)}}\big| \frac{|u|^q}{|x|^{\frac{q\eta}{p^*(\eta,\mu)}}}dx\\
& \leq  \Big(\int_\Omega \big|a(x)\big|^{\frac{p^*(\eta, \mu)}{p^*(\eta,\mu)-q}}|x|^{\frac{q\eta}{p^*(\eta,\mu)-q}}dx\Big)^{\frac{p^*(\eta,\mu)-q}{p^*(\eta,\mu)}}
\Big(\int_\Omega \frac{|u|^{p^*(\eta,\mu)}}{|x|^\eta}dx\Big)^{\frac{q}{p^*(\eta,\mu)}}.
\end{align*}
By the Hardy-Sobolev inequality (\ref{2013-9-05-e5}) and the assumption that $\big|a(x)\big|^{\frac{p^*(\eta, \mu)}{p^*(\eta,\mu)-q}}|x|^{\frac{q\eta}{p^*(\eta,\mu)-q}}$ $ \in L^1(\Omega)$, we see  that $\chi(u)$ is well defined.  Now we assume that $u_n\rightharpoonup u$ in $W_{\mu}^{1,p}(\Omega)$. By Lemma     \ref{2014-9-6-l1},  going to a subsequence if necessary,   we may assume that
$$u_n\rightarrow u\;a.e.\;\hbox{on}\;\Omega.$$
By the Hardy-Sobolev inequality again, we see that $\{u_n\}$ is bounded in $L^{p^*(\eta,\mu)}(\Omega, \frac{dx}{|x|^\eta})$, then $\{\frac{|u_n|^q}{|x|^{\frac{q\eta}{p^*(\eta,\mu)}}}\}$ is bounded in $L^{\frac{p^*(\eta,\mu)}{q}}(\Omega)$. Hence $\frac{|u_n|^q}{|x|^{\frac{q\eta}{p^*(\eta,\mu)}}}\rightharpoonup \frac{|u|^q}{|x|^{\frac{q\eta}{p^*(\eta,\mu)}}}$ in $L^{\frac{p^*(\eta,\mu)}{q}}(\Omega)$ up to a subsequence. Recalling that $a(x)|x|^{\frac{q\eta}{p^*(\eta,\mu)}}\in L^{\frac{p^*(\eta,\mu)}{p^*(\eta,\mu)-q}}$ and
$$\frac{1}{\frac{p^*(\eta,\mu)}{q}}+\frac{1}{\frac{p^*(\eta,\mu)}{p^*(\eta,\mu)-q}}=1, $$
we obtain that
$$\int_\Omega a(x)|u_n|^pdx\rightarrow \int_\Omega a(x)|u|^pdx.$$
Thus, we prove that $\chi$ is weakly continuous. Especially, when $\mu\leq 0<\mu+p<N$ and $q=p$ , we can take $\eta=0$ and obtain the final result.
\ep
\br\lab{2014-10-16-r1}
 When $q=p=2, \mu=0$,    Lemma \ref{2014-9-8-l1}  is exactly the Lemma 2.13 in  \cite{Willem.1996}. Evidently, such a very typical case is essentially different from  the general situation considered here.
\er

\vskip 0.3in
Consider the  minimizing problem
\begin{itemize}
\item[$(Q)$]  $\min \displaystyle \Big\{\int_\Omega \frac{|\nabla u|^p}{|x|^\mu}dx: \;  u\in W_{\mu}^{1,p}(\Omega), \int_\Omega V|u|^qdx=1\Big\}$.
\end{itemize}


\bl\lab{2014-9-8-l2}
Under  the assumption ${\bf (H)}$, $\int_\Omega V_+|u|^qdx$ is weakly continuous in $W_{\mu}^{1,p}(\Omega)$.
\el
\bp The proof is inspired by that of  \cite[Lemma 2.1]{SzulkinWillem.1999}. However, our case is much more complicated. In view of Lemma \ref{2014-9-8-l1}, we only need to prove that $\int_\Omega V_2|u|^pdx$ is weakly continuous.

\vskip0.1in

\noindent {\bf Step 1. } We prove that $\{V_2|u_n|^q\}$ is a tight sequence. We only need to prove the case  of  $(H_2)$ or $(H_3)$.
For the case of $(H_2)$, by the H\"{o}lder inequality and the Hardy-Sobolev inequality (\ref{2013-9-05-e5}), we see that
\begin{align*}
&\int_{\Omega\cap B_R^c(0)} V_2|u_n|^qdx=\int_{\Omega\cap B_R^c(0)} V_2|x|^{\frac{(\mu+p)q}{p}} \frac{|u_n|^q}{|x|^{\frac{(\mu+p)q}{p}}}dx\\
& \leq  \Big(\int_{\Omega\cap B_R^c(0)} \big(V_2(x)\big)^{\frac{p}{p-q}}|x|^{\frac{(\mu+p)q}{p-q}}dx\Big)^{\frac{p-q}{p}} \Big(\int_{\Omega\cap B_R^c(0)}\frac{|u_n|^p}{|x|^{\mu+p}}dx
\Big)^{\frac{q}{p}}\\
& \rightarrow  0\;\hbox{as}\;R\rightarrow \infty.
\end{align*}
For the case of $(H_3)$, since $q\geq p, \mu+p<N$, we have $\bar{\sigma}:=N-\frac{q}{p}(N-\mu-p)\leq \mu+p$ and $q=p^*(\bar{\sigma},\mu)$.
Recall that  $\displaystyle \lim_{\stackrel{|x|\rightarrow \infty}{x\in\Omega}}|x|^{\bar{\sigma}}V_2(x)=0$ and the Hardy-Sobolev inequality, we see that
$\{V_2|u_n|^q\}$ is also a tight sequence.
In summary, under the assumption (H), for $\forall\;\varepsilon>0$, we can take $R>0$ large enough such that
\be\lab{2014-9-8-e2}
\int_{\Omega\backslash B_R(0)} V_2|u_n|^qdx<\varepsilon\;\hbox{for all}\;n.
\ee
It follows from Fatou's Lemma, we also have
\be\lab{2014-9-8-e3}
\int_{\Omega\backslash B_R(0)} V_2|u|^qdx<\varepsilon.
\ee

\noindent {\bf Step 2. }
We note that for the cases of $(H_1)$ and $(H_2)$, $1<q<p$.
Due to the compactness, we can choose a finite covering of $\overline{\Omega\cap B_R(0)}$ by closed balls $\overline{B_{r_i}(x_i)}, 1\leq i\leq k$ such that $\{x_i\}_{i=1}^{k}\subset \overline{\Omega\cap B_R(0)}\subset \bar{\Omega}$ and that
\be\lab{2014-9-8-e4}
|x-x_i|^{\mu+p}V_2(x)\leq C\;\hbox{for all}\;x\in B_{r_i}(x_i), i=1,2,\cdots,k.
\ee
Note  that $1<q<p=p^*(\mu+p, \mu)$, then by the weighted Rellich-Kondrachov compactness Theorem \ref{2014-9-7-th1}, it is easy to obtain that
\be\lab{2014-9-10-xe1}
\int_{\Omega\cap B_R(0)}V_2|u_n|^qdx\rightarrow \int_{\Omega\cap B_R(0)}V_2|u|^qdx\;\hbox{as}\;n\rightarrow \infty.
\ee
Hence, for the cases of $(H_1)$ and $(H_2)$, by (\ref{2014-9-8-e2}), (\ref{2014-9-8-e4}) and (\ref{2014-9-10-xe1}), we prove that
$$\int_\Omega V_2|u_n|^qdx\rightarrow \int_\Omega V_2|u|^qdx\;\hbox{as}\;n\rightarrow \infty.
$$
For the case of $(H_3)$, since $q\geq p, \mu+p<N$, we see that $p^*(\bar{s}, \mu)>\max\{p,q\}$. By compactness again, for $\forall\;\varepsilon>0$,  we can choose a finite covering of $\overline{\Omega\cap B_R(0)}$ by closed balls $\overline{B_{r_i}(x_i)}, 1\leq i\leq k$ such that $\{x_i\}_{i=1}^{k}\subset \overline{\Omega\cap B_R(0)}\subset \bar{\Omega}$ and
\be\lab{2014-9-10-xe2}
|x-x_i|^{\bar{s}}V_2(x)\leq \varepsilon\;\hbox{for all}\;x\in B_{r_i}(x_i),\;\;  i=1,2,\cdots,k.
\ee
We note that $k$ depends on $\varepsilon$.
 By the assumption $(H_{3})$ again, we can take $0<r<\min\{r_1,r_2,\cdots, r_k\}$ such that
\be\lab{2014-9-8-e5}
|x-x_i|^{\bar{s}}V_2(x)\leq \frac{\varepsilon}{k}\;\hbox{for all}\;x\in B_{r}(x_i),\;\;  i=1,2,\cdots,k.
\ee
Set $$A:=\bigcup_{i=1}^{k}B_r(x_i),$$
then by the Hardy-Sobolev inequality (\ref{2013-9-05-e5}),
\be\lab{2014-9-10-xe3}
\Big(\int_{B_r(x_i)} \frac{|u_n|^q}{|x-x_i|^{\bar{s}}}dx\Big)^{\frac{1}{q}}\leq C\Big(\int_{B_r(x_i)} \frac{|\nabla u_n|^p}{|x-x_i|^{\mu}}dx\Big)^{\frac{1}{p}},\;i=1,2,\cdots,k.
\ee
Thus,
\be\lab{2014-9-8-e6}
\int_A V_2|u_n|^qdx\leq \varepsilon C^q, \int_A V_2|u|^qdx\leq \varepsilon C^q.
\ee
It follows from (\ref{2014-9-10-xe2}) that $V_2\in L^\infty\big((\Omega\cap B_R(0))\backslash A\big)$, and then $V_2$ satisfies the assumption of Lemma \ref{2014-9-8-l1} up to the bounded domain $(\Omega\cap B_R(0))\backslash A$. So
\be\lab{2014-9-8-e7}
\int_{(\Omega\cap B_R(0))\backslash A} V_2|u_n|^qdx\rightarrow \int_{(\Omega\cap B_R(0))\backslash A} V_2|u|^qdx.
\ee
Then, by (\ref{2014-9-8-e2}), (\ref{2014-9-8-e3}), (\ref{2014-9-8-e6}) and (\ref{2014-9-8-e7}), we also have
$$\int_\Omega V_2|u_n|^qdx\rightarrow \int_\Omega V_2|u|^qdx.$$
\ep
\bc\lab{2014-9-8-cor1}
Under the assumption $(H)$, $\int_\Omega V(x)|u|^qdx$ is weakly  upper semicontinuous in $W_{\mu}^{1,p}(\Omega)$.
\ec
\bp
It is an obvious  conclusion which can be deduced  by
Lemma \ref{2014-9-8-l2} and the Fatou's Lemma.
\ep

\bt\lab{2014-9-8-th1}
Under the assumption $(H)$, problem $(Q)$ has a solution $\varphi_1\geq 0$. Moreover, $(\varphi_1,\lambda_1)$ is a solution to problem (\ref{2014-9-7-xe2}), where  $\lambda_1:=\int_\Omega \frac{|\nabla \varphi_1|^p}{|x|^\mu}dx$.
\et
\bp
Let $\{u_n\}$ be a minimizing sequence for $(Q)$. By Lemma \ref{2014-9-6-l1}, we may assume that $u_n\rightharpoonup u$ in $W_{\mu}^{1,p}(\Omega)$ and $u_n\rightarrow u$ a. e.  on $\Omega$. Hence,
$$\int_\Omega\frac{|\nabla u|^p}{|x|^\mu}dx\leq \liminf_{n\rightarrow \infty}\int_\Omega \frac{|\nabla u_n|^p}{|x|^\mu}dx=\inf(Q).$$
By Corollary \ref{2014-9-8-cor1}, we have that $\displaystyle \int_\Omega V|u|^qdx\geq 1.$
Let
$$\varphi_1:=\frac{u}{(\int_\Omega V|u|^qdx)^{\frac{1}{q}}},$$
we see that
$\int_\Omega V|\varphi_1|^qdx=1$ and
$$\inf(Q)\leq\int_\Omega \frac{|\nabla \varphi_1|^p}{|x|^\mu}dx=\frac{1}{(\int_\Omega V|u|^qdx)^{\frac{p}{q}}}\int_\Omega \frac{|\nabla u|^p}{|x|^\mu}dx\leq \int_\Omega \frac{|\nabla u|^p}{|x|^\mu}dx\leq \inf(Q).$$
Hence, we see that $\int_\Omega V|u|^qdx=1$ and $\varphi_1=u$
is a solution of $(Q)$. Note that $|\varphi_1|$ is also a solution, we may assume $\varphi_1\geq 0$.  Moreover,  there exists some Lagrange multiplier $\lambda_1$ such that
$$- div(\frac{|\nabla \varphi_1|^{p-2}\nabla \varphi_1}{|x|^\mu})=\lambda_1 V(x)|\varphi_1|^{q-2}\varphi_1.$$
Testing by $\varphi_1$, we have
$$\int_\Omega \frac{|\nabla \varphi_1|^p}{|x|^\mu}dx=\lambda_1 \int_\Omega V(x)|\varphi_1|^qdx=\lambda_1>0.$$
We also note that $\lambda_1=\inf(Q).$
\ep

We need the following Br\'{e}zis-Lieb type lemma.

\bl\lab{2014-9-8-l3}
Let $\Omega$ be an open subset of  $ \R^N$ and
assume that $\{u_n\}$ satisfies   $$q\geq 1, \sup_{n}\int_\Omega |a(x)||u_n|^q dx<\infty\;\hbox{and}\;u_n\rightarrow u\; a.e. \;\hbox{in}\;\Omega.$$
 Then
\be\lab{zz=2}\lim_{n\rightarrow \infty} \int_\Omega a(x)(|u_n|^q-|u_n-u|^q)=\int_\Omega a(x)|u|^qdx.\ee
\el
\bp
Consider the new measure $\nu$ such that $d\nu=|a(x)|dx$, then    (\ref{zz=2})    can be deduced from  the Br\'{e}zis-Lieb Lemma  with respect   to the new measure $\nu$:
$$\lim_{n\rightarrow \infty}\int_\Omega\Big||u_n|^q-|u_n-u|^q-|u|^q\Big|d\nu=0.$$
\ep

Now we introduce a new space $E:=\{u\in W_{\mu}^{1,p}(\Omega): \|u\|_E<\infty\}$, where
$$\|u\|_E:=\Big(\int_\Omega \frac{|\nabla u|^p}{|x|^\mu}dx\Big)^{\frac{1}{p}}+\Big(\int_\Omega V_- |u|^qdx\Big)^{\frac{1}{q}}.$$
Set
$$I(u):=\int_\Omega \frac{|\nabla u|^p}{|x|^\mu}dx,\quad  J(u):=\int_\Omega V|u|^qdx,\quad J_\pm(u):=\int_\Omega V_\pm |u|^qdx,$$
then we see that
$M:=\{u\in E:   J(u)=1\}$
is a $C^1$-manifold.

\bl\lab{2014-9-8-l4}
If $V$ satisfies $(H)$, then
\begin{itemize}
\item[$(i)$] there exists some $C>0$ such that $J_+(u)\leq C\big(I(u)\big)^{\frac{q}{p}}$ for all $u\in E$;
\item[$(ii)$]  $J_+$ is weakly continuous and $J'_+$ is completely continuous (or weak-to-strong continuous), i.e., if $u_n\rightharpoonup u$, then $J'_+(u_n)\rightarrow J'_+(u)$.
\end{itemize}
\el
\bp
Let $u_n\rightharpoonup u$ in $E$, since $E\subset W_{\mu}^{1,p}(\Omega)$, we may assume that $u_n\rightarrow u$ a.e. in $\Omega$. We note that $(i)$ is deduced by Lemma  \ref{2014-9-8-l2}.
Next we shall prove $(ii)$. By Lemma \ref{2014-9-8-l3}, we see that $J_+$ is weakly continuous.
For any $v\in E$, by   the H\"{o}lder inequality  up  to the new measure $d\nu=V_+dx$,   we have
\begin{align}\lab{2014-9-8-e9}
\Big|\big\langle & J'_+(u_n)-J'_+(u), v\big\rangle\Big| \nonumber \\=&\Big|\int_\Omega V_+\big(|u_n|^{q-2}u_n-|u|^{q-2}u\big)vdx\Big|\nonumber\\
=&\Big|\int_\Omega \big(|u_n|^{q-2}u_n-|u|^{q-2}u\big)vd\nu\Big|\nonumber\\
\leq&\Big(\int_\Omega \big||u_n|^{q-2}u_n-|u|^{q-2}u\big|^{\frac{q}{q-1}}d\nu\Big)^{\frac{q-1}{q}} \Big(\int_\Omega |v|^qd\nu\Big)^{\frac{1}{q}}\nonumber\\
\leq&C\|v\|_E \Big(\int_\Omega \big||u_n|^{q-2}u_n-|u|^{q-2}u\big|^{\frac{q}{q-1}}d\nu\Big)^{\frac{q-1}{q}}.
\end{align}
Let $v_n:=|u_n|^{q-2}u_n, v:=|u|^{q-2}u$. Since $\frac{q}{q-1}>1$,   by Lemma \ref{2014-9-8-l3}    we have
\begin{align}\lab{2014-9-8-e10}
&\lim_{n\rightarrow \infty}\int_\Omega \big||u_n|^{q-2}u_n-|u|^{q-2}u\big|^{\frac{q}{q-1}}d\nu\nonumber\\
&= \lim_{n\rightarrow \infty} \int_\Omega |v_n-v|^{\frac{q}{q-1}}d\nu\nonumber\\
& =\lim_{n\rightarrow \infty}\int_\Omega \big(|v_n|^{\frac{q}{q-1}}-|v|^{\frac{q}{q-1}}\big)d\nu\nonumber\\
& = \lim_{n\rightarrow \infty} \int_\Omega \big(|u_n|^q-|u|^q\big)d\nu\nonumber\\
& = 0\;\hbox{since}\;J_+\;\hbox{is weakly continuous}.
\end{align}
It follows  from (\ref{2014-9-8-e9}) and (\ref{2014-9-8-e10}) that
$J'_+$ is completely continuous.
\ep

\vskip0.2in

\noindent{\bf Proof of Theorem  \ref{2014-9-7-th3}. }   For the convenience of the readers, we recall the  Krasnoselskii Genus. Define  $\mathcal{A}=\{A\subset M:   A\;\hbox{closed}, A=-A\}$.
For $A\subset \mathcal{A}, A\neq \emptyset$, let
$$\gamma(A):=\begin{cases}\inf\{m: \; \;\exists\;h\in C^0(A;\R^m\backslash \{0\}), h(-u)=-h(u)\}, \\ \infty,\quad \hbox{if}\;\{\cdots\}=\emptyset,\;\hbox{in particular, if}\; 0\in A, \end{cases}$$
and let $\gamma(\emptyset)=0$.
Define
$$\lambda_n:=\inf_{\gamma(A)\geq n}\sup_{u\in A}I(u), \quad n=1,2,\cdots.$$
Under the assumption $(H)$, we see that $\{x\in \Omega: V(x)>0\}$ has positive measure.  Note that  $\gamma(\mathcal{S}^{n-1})=n$, where $\mathcal{S}^{n-1}$ is the unit sphere of $\R^n$, it follows that $\lambda_n$ is well defined for all $n$ by constructing a suitable odd homeomorphism. Moreover, we see that $\lambda_1=\inf_{u\in M}I(u)>0$ coincides with the value given by Theorem \ref{2014-9-8-th1}.   We now prove that  $I\Big|_M$ satisfies $PS$ condition. Let $\{u_n\}\subset E $ be a $PS$ sequence. Then there is a corresponding   sequence  $\mu_k\in \R$ such that
\be\lab{2014-9-8-xe1}
A_{\mu_k}(u_k):=I'(u_k)-\mu_k J'(u_k)\rightarrow 0\;\hbox{in}\;E^*.
\ee
It follows  from  $(i)$ of Lemma \ref{2014-9-8-l4} that $J_+(u_k)$ is bounded and therefore $J_-(u_k)$ is bounded  since $
J_-(u_k)=J_+(u_k)-1.
$
Note that $\|u_k\|_E\equiv I(u_k)^{\frac{1}{p}}+J_-(u_k)^{\frac{1}{q}}$, we see that $\{u_k\}$ is bounded in $E$ and then $J'(u_k)$ is bounded. Up to a subsequence, we may assume that $u_k\rightharpoonup u$ in $E$ and $u_k\rightarrow u$ a.e.  in $\Omega$.
By Corollary \ref{2014-9-8-cor1}, we have
$J(u)\geq 1$
and it follows that
\be\lab{2014-9-8-xe3}
 \langle J'(u), u\rangle=q J(u)\geq q.
\ee
Testing by $u_k$ in (\ref{2014-9-8-xe1}),  we get that
\be\lab{2014-9-8-xe4}
\langle I'(u_k), u_k\rangle -\mu_k \langle J'(u_k), u_k\rangle
=p I(u_k)-q\mu_k \rightarrow 0.
\ee
By the boundedness of $I(u_k)$ and (\ref{2014-9-8-xe4}), we obtain that $\{\mu_k\}$ is bounded.
Up to a subsequence, we assume that $\mu_k\rightarrow \mu_\infty$ and it follows from (\ref{2014-9-8-xe4}) again, we have $I(u_k)\rightarrow \frac{q}{p}\mu_\infty$ and
$A_{\mu_\infty}(u_k)\rightarrow 0, \langle A_{\mu_\infty}(u_k), u_k-u\rangle\rightarrow 0$.
By the Br\'{e}zis-Lieb Lemma, it is easy to see that
\be\lab{2014-9-8-xe5}
I(u_k-u)=I(u_k)-I(u)+o(1).
\ee
Insert  Lemma \ref{2014-9-8-l3}, we have
\be\lab{2014-9-8-xe6}
J_-(u_k-u)=J_-(u_k)-J_-(u)+o(1).
\ee
A direct calculation shows that
$$\langle A_{\mu_\infty}(u_k), u_k-u\rangle=p\Big[I(u_k)-I(u)\Big]+q\mu_\infty\Big[\big(J_-(u_k)-J_-(u)\big)\Big]. $$
Combine with  (\ref{2014-9-8-xe5}) and (\ref{2014-9-8-xe6}) that
\be\lab{2014-10-16-xe1}
pI(u_k-u)+q\mu_\infty J_-(u_k-u)\rightarrow 0.
\ee
By the weakly lower semicontinuity of a norm and (\ref{2014-9-8-xe4}), we have
$$I(u)\leq \liminf_{k\rightarrow \infty}I(u_k)=\frac{q}{p}\mu_\infty.$$
Since $J(u)\geq 1$, we have $u\neq 0$ and then $I(u)>0$.
Thus, $\mu_\infty>0$ and   further by \eqref{2014-10-16-xe1}, we have $I(u_k-u)\rightarrow 0, J_-(u_k-u)\rightarrow 0$.
Hence,
$$\|u_k-u\|_E\equiv I(u_k-u)^{\frac{1}{p}}+J_-(u_k-u)^{\frac{1}{q}}\rightarrow 0.$$
We have thus proved that $I\Big|_M$ satisfies $PS$ condition.
Hence, $\lambda_n s$ are critical values due to  \cite[page 98, Theorem 5.7]{Struwe.2008}. There exists a critical point $\varphi_n$ such that $I'(\varphi_n)=\widetilde{\mu} J'(\varphi_n)$. Testing by $\varphi_n$, we have
$$pI(\varphi_n)=\langle I'(\varphi_n), \varphi_n\rangle =\widetilde{\mu}\langle J'(\varphi_n), \varphi_n\rangle=q\widetilde{\mu} J(\varphi_n)=q\widetilde{\mu}.$$
Hence,  $\lambda_n=\frac{q}{p}\widetilde{\mu}=I(\varphi_n)$. Finally, we shall prove that $\lambda_n\rightarrow \infty$.
Note that  if $\lambda_n=\lambda_{n+1}=\cdots=\lambda_{n+k-1}=\lambda$ for some $n,k$,  then the set of critical points corresponding to $\lambda$ has genus $\geq k$ (see \cite[page 97, Lemma 5.6]{Struwe.2008}).
We also note that  $\lambda_{n+1}\geq \lambda_n$,    apply the similar argument of \cite[Proposition 9.33]{Rabinowitz.1986}, we can prove that $\lambda_n\rightarrow \infty$ as $n\rightarrow \infty$.  Moreover, if $q=p$, we see that $\lambda_n s$ are eigenvalues such that $\lambda_n\rightarrow \infty$ as $n\rightarrow \infty$ and  $\varphi_n$ is an eigenfunction of (\ref{2014-9-7-xe2}) corresponding to $\lambda_n$.
On the other hand,  if $1<q\neq p$, after scaling we see that $v_n:=(\frac{\lambda}{\lambda_n})^{\frac{1}{p-q}}\varphi_n$ is a sequence of solutions to  \eqref{2014-9-7-xe2} such that
\begin{align}
c_n:=&\Phi(v_n)=:\frac{1}{p}I(v_n)-\frac{1}{q}\lambda J(v_n)\nonumber\\
=&(\frac{1}{p}-\frac{1}{q})I(v_n)\nonumber\\
=&(\frac{1}{p}-\frac{1}{q})(\frac{\lambda}{\lambda_n})^{\frac{p}{p-q}}I(u_n)\nonumber\\
=&(\frac{1}{p}-\frac{1}{q})(\frac{\lambda}{\lambda_n})^{\frac{p}{p-q}}\lambda_n\nonumber\\
=&(\frac{1}{p}-\frac{1}{q})\lambda^{\frac{p}{p-q}} \lambda_{n}^{\frac{q}{q-p}}.
\end{align}
Hence, if $q>p$, \eqref{2014-9-7-xe2} possesses a sequence of solutions with energy $0<c_1\leq c_2\leq \cdots\leq c_n\rightarrow \infty$ as $n\rightarrow \infty$. If $q<p$, \eqref{2014-9-7-xe2} has  a sequence of solutions with energy $-\infty<c_1\leq c_2\leq \cdots\leq c_n\rightarrow 0$ as $n\rightarrow \infty$.
\hfill$\Box$

\vskip 0.3in
As an application, we obtain the following result:
\bc\lab{2014-9-10-cro1}
Assume that $1<p<N, \mu+p<N, 1<q<p^*$ and
let $\Omega\subset \R^N$ be an open Lipschitz domain such that one of the following holds:
\begin{itemize}
\item[$(i)$] $\Omega$ is bounded and $1<q<p^*(\sigma, \mu)$.
\item[$(ii)$] $\Omega$ has finite Lebergue measure, $0\leq \sigma$ and $1<q<p^*(\sigma,\mu)$.
\end{itemize}
Consider the following problem
\be\lab{2014-9-10-we1}
-div\big(\frac{|\nabla u|^{p-2}\nabla u}{|x|^\mu}\big)=\lambda \frac{1}{|x|^\sigma}|u|^{q-2}u\;\hbox{ in }\; \Omega,\;
u\in W_{\mu}^{1,p}(\Omega).
\ee

 \begin{itemize}

\item[(1)] If $q=p$, then  (\ref{2014-9-10-we1})
 possesses a sequence of eigenvalues $0<\lambda_1\leq \lambda_2\leq \cdots\leq \lambda_n\rightarrow \infty$.

 \item[(2)] If $q>p$,   then   for any positive fixed $\lambda$, \eqref{2014-9-10-we1} possess a sequence of solutions $\{v_n\}$ such that $0<c_1\leq c_2\leq \cdots\leq c_n\rightarrow \infty$ as $n\rightarrow \infty$.

\item[(3)] If $q<p$,  then  for any positive fixed $\lambda$, \eqref{2014-9-10-we1} possess a sequence of solutions $\{v_n\}$ such that  $-\infty<c_1\leq c_2\leq\cdots\leq c_n\rightarrow 0$ as $n\rightarrow \infty$,
 \end{itemize}
where $$c_n:=\frac{1}{p}\int_\Omega \frac{|\nabla v_n|^p}{|x|^\mu}dx-\frac{1}{q}\lambda\int_\Omega\frac{|v_n|^q}{|x|^\sigma}dx.$$
\ec
\bp
By Theorem\ref{2014-9-7-th3}, we only need to check  that $V(x)=|x|^{-\sigma}$ satisfies the assumption $(H)$. We first check the case $(i)$:
when $\Omega$ is bounded and $p<p^*(\sigma, \mu)$, we have $\sigma<\mu+p<N$, then we take $V(x)\equiv V_1(x)+V_2=|x|^{-\sigma}+0$ and choose  $\eta=\sigma$ in the assumption $(H)$, we see that
$\frac{-p^*(\sigma, \mu)\sigma+q\sigma}{p^*(\sigma,\mu)-q}=-\sigma>-N$.
Thus,
$\big|V(x)\big|^{\frac{p^*(\eta, \mu)}{p^*(\eta,\mu)-q}}|x|^{\frac{q\eta}{p^*(\eta,\mu)-q}}\in L^1(\Omega)$.
When $\Omega$ is bounded and $p^*(\sigma, \mu)\leq p$, we have that $\mu+p\leq \sigma$. By $q<p^*(\sigma, \mu)$, we have $\sigma<N-\frac{q}{p}(N-\mu-p)<N$. Hence, we take $V(x)\equiv V_1(x)+V_2(x)=|x|^{-\sigma}+0$,
recalling that $q<p^*(\sigma,\mu)$, we have
\be\lab{2014-9-10-wwe1}
\frac{-p^*(\eta, \mu)\sigma+q\eta}{p^*(\eta,\mu)-q}>-N\Leftrightarrow \eta<N.
\ee
Then for any $\eta$ satisfying  $\frac{N\mu}{N-p}\leq \eta<\mu+p$, $\big|V(x)\big|^{\frac{p^*(\eta, \mu)}{p^*(\eta,\mu)-q}}|x|^{\frac{q\eta}{p^*(\eta,\mu)-q}}\in L^1(\Omega)$. The assumption $(H)$ holds.

Next, we check for the case $(ii)$: we prefer to introduce the characteristic function
for any subset  $ A\subset \R^N:$
$$1_A(x)=\begin{cases}1,\;&x\in A,\\ 0,&x\not\in A. \end{cases}\;$$
{\bf Case 1. }If $p^*(\sigma,\mu)\leq p$, we take $V(x)\equiv V_1(x)+V_2(x)$ with
$$V_1(x)=|x|^{-\sigma} 1_{\Omega\cap B_r(0)}(x), \quad  V_2(x)=|x|^{-\sigma} 1_{\Omega\cap B_r^c(0)}(x)$$
for some $r>0$. Then similar to the arguments above, we have that $$\big|V_1(x)\big|^{\frac{p^*(\eta, \mu)}{p^*(\eta,\mu)-q}}|x|^{\frac{q\eta}{p^*(\eta,\mu)-q}}\in L^1(\Omega), $$
since $q<\min\{p^*(\sigma, \mu),p^*\}$.
 Note that in this case, we have $\sigma\geq \mu+p$ and $p>q>1$.
It follows  from $\sigma\geq 0$ that
$-\sigma p+(\mu+p)q=-\sigma(p-q)-(\sigma-\mu-p)q\leq 0$. Hence
\be\lab{2014-9-10-wwwe1}
\lim_{R\rightarrow \infty}\int_{\{x\in \Omega:|x|>R\}} \big(V_2(x)\big)^{\frac{p}{p-q}}|x|^{\frac{(\mu+p)q}{p-q}}dx=0
\ee
as $R\rightarrow \infty$ due to the fact $L(\Omega)<\infty$.
Hence,  the assumption $(H)$ is satisfied.

\vskip0.1in

\noindent {\bf Case 2. } If $p^*(\sigma,\mu)> p$, we have $0\leq \sigma<\mu+p<N$ and $q<p^*(\sigma, \mu)\leq p^*(0,\mu)$. Now, take $V(x)\equiv V_1(x)+V_2(x)$ with
$$V_1(x)=|x|^{-\sigma} 1_{\Omega\cap B_r^c(0)}(x),\quad  V_2(x)=|x|^{-\sigma} 1_{\Omega\cap B_r(0)}(x)$$
for some $r>0$.
For $V_1(x)$, take $\eta=\sigma$, we have
$$\big|V_1(x)\big|^{\frac{p^*(\eta, \mu)}{p^*(\eta,\mu)-q}}|x|^{\frac{q\eta}{p^*(\eta,\mu)-q}}=1_{\Omega\cap B_r^c(0)}(x)|x|^{-\sigma}\in L^1(\Omega)$$ since $\sigma\geq 0$ and $L(\Omega)<\infty$. We also note that
\begin{itemize}
\item[(1)]if $q<p$, $V_2$ satisfies $(H_2)$ since $\sigma<\mu+p$.
\item[(2)]if $p\leq q<p^*(\sigma,\mu)$, we have $\bar{\sigma}:=N-\frac{q}{p}(N-\mu-p)>\sigma$. Hence, $V_2$ satisfies $(H_3)$.
\end{itemize}
Hence, the assumption $(H)$ is also satisfied for the case of $(ii)$.
\ep


\vskip 0.3in
\noindent {\bf Acknowledgment.} The authors would like to thank professor G. Cerami  for her encouragement, discussions and pointing out some references in the preparation of the manuscript.


\vskip0.26in

\end{CJK*}
 \end{document}